\documentclass[12pt]{article}

\newcommand{\qed}{{\hfill$\square$}\medskip}

\usepackage{amsmath, epsfig, cite}
\usepackage{amssymb}
\usepackage{amsfonts}
\usepackage{latexsym}

\newtheorem{thm}{Theorem}[section]

\newtheorem{cor}[thm]{Corollary}

\begin{document}

\begin{center}
{\Large\bf A recursive algorithm for trees and forests}
\end{center}

\vskip 2mm \centerline{Song Guo and Victor J. W. Guo\footnote{Corresponding author.}}
\begin{center}
{School of Mathematical Sciences, Huaiyin Normal University, Huai'an, Jiangsu 223300,
 People's Republic of China\\[5pt]
{\tt guosong77@hytc.edu.cn,  jwguo@hytc.edu.cn} }

\end{center}


\vskip 0.7cm

\noindent{\bf Abstract.} Trees or rooted trees have been generously studied in the literature.
A forest is a set of trees or rooted trees. Here we give recurrence relations between the
number of some kind of rooted forest with $k$ roots and that with $k+1$ roots on
$\{1,2,\ldots,n\}$. Classical formulas for counting  various  trees such as
rooted trees, bipartite trees, tripartite trees, plane trees, $k$-ary plane
trees, $k$-edge colored trees follow immediately from our recursive relations.


\vskip 2mm
\noindent{\it MR Subject Classifications}: 05C05, 05A15, 05A19.


\section{Introduction and notations}
The famous Cayley's formula for counting trees states that the number of labeled trees
on $[n]:=\{1,\ldots,n\}$ is $n^{n-2}$. Clarke \cite{Clarke} first gives a refined version for Cayley's
formula by setting up a recurrence relation. Erd$\acute{\rm e}$lyi and Etherington
\cite{EE} gives a bijection between semilabeled trees and partitions. This bijection was
also discovered by Haiman and Schmitt \cite{HS}. A general bijective algorithm was
given by Chen \cite{Chen90}. Aigner and Ziegler's book \cite{AZ} collected four different proofs
of Cayley's formula. We refer the reader to \cite{ChenGuo,ChenPeng,DuYin,Hou,GZ,SZ} for
several recent results on the enumeration of trees. The goal of this paper is to establish simple linear recurrences
between certain forests with roots $1,\ldots,k$ and forests with roots $1,\ldots,k+1$,
from which one can deduce several classical results on counting trees.

The set of forests of $k$ rooted trees on $[n]$ with roots $1,\ldots,k$ is denoted by
$\mathcal{F}_n^k$. Suppose that $F\in\mathcal{F}_n^k$ and $x$ is a vertex of $F$, the subtree
rooted  at $x$ is denoted by $F_x$. We say that a vertex $y$ of $F$ is a {\it descendant}
of $x$, if $y$ is a vertex of $F_x$, i.e., $x$ is on the path from the root of $T$ to the
vertex $y$, and is denoted by $y\prec x$.

For any edge $e=(x,y)$ of a tree $T$ in a forest $F$, if $y$ is a vertex of $T_x$, we call
$x$ the {\it father vertex} of $e$, $y$ the {\it child vertex} of $e$, $x$ the father of $y$,
and $y$ a child of $x$, sometimes we also say $e$ is {\it out of} $x$.
The degree of a vertex $x$ in a rooted tree $T$ is the number
of children of $x$, and is denoted by $\deg_{T}(x)$, or $\deg_{F}(x)$.
As usual, a vertex with degree zero is called a leaf.

An unrooted labeled tree will be treated as a rooted tree in which the smallest vertex is
chosen as the root. Moreover, if $\mathcal{A}$ is a set of trees, then we will use
$\mathcal{A}[P]$ to denote the subset of all elements of  $\mathcal{A}$ satisfying the condition $P$.
\section{The fundamental recursion}
One of our main results is as follows.
\begin{thm}\label{thm-main}
For $2\leqslant k\leqslant n-1$, we have the following recurrence relation:
\begin{equation}\label{eq-main}
|\mathcal{F}_n^{k-1}[n\prec 1]|=n|\mathcal{F}_n^{k} [n\prec 1]|.
\end{equation}
\end{thm}

\noindent{\it Proof.} Suppose that $F\in\mathcal{F}_n^{k-1}[n\prec 1]$. First,
remove the subtree $F_k$ from $F$ and add it to be a new tree
in the forest. Second if $n$ is not a descendant of $1$ in the new
forest, $n$ must be a descendant of $k$, exchange labels of the vertices $1$
and $k$. Thus, we obtain a forest $F'\in\mathcal{F}_n^k[n\prec
1]$.

Conversely, for a forest $F'\in\mathcal{F}_n^k[n\prec 1]$, we can
attach $F'_k$ to any vertex of the other trees in
$F$ as a subtree, or attach $F'_1$ to any vertex of
$F'_k$ as a subtree and exchange labels of the vertices $1$ and $k$.
The proof then follows from the fact that $F'$ has $n$ vertices
altogether. \qed \vskip 0.5cm
\begin{figure}[h]
\begin{center}
\setlength{\unitlength}{0.08in}
\begin{tabular}{|c|}\hline
\begin{picture}(10,18)
\put(4,14){\line(-1,-2){2}}
\put(4,14){\line(1,-2){2}}
\put(6,10){\line(0,-1){4}}
\put(4,14){\circle*{0.7}}        \put(3.5,15){\makebox(2,1)[l]{1}}
\put(2,10){\circle*{0.7}}        \put(1.6,8){\makebox(2,1)[l]{5}}
\put(6,10){\circle*{0.7}}        \put(6.7,9.5){\makebox(2,1)[l]{3}}
\put(6,6){\circle*{0.7}}         \put(5.6,4){\makebox(2,1)[l]{4}}
\put(8,14){\circle*{0.7}}        \put(7.6,15){\makebox(2,1)[l]{2}}
\end{picture}
   \vline
\begin{picture}(10,18)
\put(3,14){\line(0,-1){12}}
\put(3,14){\circle*{0.7}}          \put(2.5,15){\makebox(2,1)[l]{1}}
\put(3,10){\circle*{0.7}}          \put(3.7,9.5){\makebox(2,1)[l]{5}}
\put(3,6){\circle*{0.7}}           \put(3.7,5.5){\makebox(2,1)[l]{3}}
\put(3,2){\circle*{0.7}}           \put(3.7,1.5){\makebox(2,1)[l]{4}}
\put(7,14){\circle*{0.7}}          \put(6.6,15){\makebox(2,1)[l]{2}}
\end{picture}
    \vline
\begin{picture}(10,18)
\put(3,14){\line(0,-1){4}}
\put(7,14){\line(0,-1){8}}
\put(3,14){\circle*{0.7}}    \put(2.5,15){\makebox(2,1)[l]{1}}
\put(3,10){\circle*{0.7}}    \put(2.6,8){\makebox(2,1)[l]{5}}
\put(7,14){\circle*{0.7}}    \put(6.6,15){\makebox(2,1)[l]{2}}
\put(7,10){\circle*{0.7}}    \put(7.7,9.5){\makebox(2,1)[l]{3}}
\put(7,6){\circle*{0.7}}     \put(6.6,4){\makebox(2,1)[l]{4}}
\end{picture}
   \vline
\begin{picture}(10,18)
\put(3,14){\line(0,-1){12}}
\put(3,14){\circle*{0.7}}     \put(2.5,15){\makebox(2,1)[l]{1}}
\put(3,10){\circle*{0.7}}     \put(3.7,9.5){\makebox(2,1)[l]{4}}
\put(3,6){\circle*{0.7}}      \put(3.7,5.5){\makebox(2,1)[l]{3}}
\put(3,2){\circle*{0.7}}      \put(3.7,1.5){\makebox(2,1)[l]{5}}
\put(7,14){\circle*{0.7}}     \put(6.6,15){\makebox(2,1)[l]{2}}
\end{picture}
   \vline
\begin{picture}(10,18)
\put(4,14){\line(-1,-2){2}}
\put(4,14){\line(1,-2){2}}
\put(6,10){\line(0,-1){4}}
\put(4,14){\circle*{0.7}}     \put(3.5,15){\makebox(2,1)[l]{1}}
\put(2,10){\circle*{0.7}}     \put(1.6,8){\makebox(2,1)[l]{4}}
\put(6,10){\circle*{0.7}}     \put(6.7,9.5){\makebox(2,1)[l]{3}}
\put(6,6){\circle*{0.7}}      \put(5.6,4){\makebox(2,1)[l]{5}}
\put(8,14){\circle*{0.7}}     \put(7.6,15){\makebox(2,1)[l]{2}}
\end{picture}
\\ \hline
\begin{picture}(13,11)
\put(3,7){\line(0,-1){4}}
\put(11,7){\line(0,-1){4}}
\put(3,7){\circle*{0.7}}    \put(2.5,8){\makebox(2,1)[l]{1}}
\put(3,3){\circle*{0.7}}    \put(2.6,1){\makebox(2,1)[l]{5}}
\put(7,7){\circle*{0.7}}    \put(6.6,8){\makebox(2,1)[l]{2}}
\put(11,7){\circle*{0.7}}   \put(10.6,8){\makebox(2,1)[l]{3}}
\put(11,3){\circle*{0.7}}   \put(10.6,1){\makebox(2,1)[l]{4}}
\end{picture}
\\ \hline
\end{tabular}
\caption{An example of Theorem \ref{thm-main} for $n=5$ and $k=3$.}
\end{center}
\end{figure}
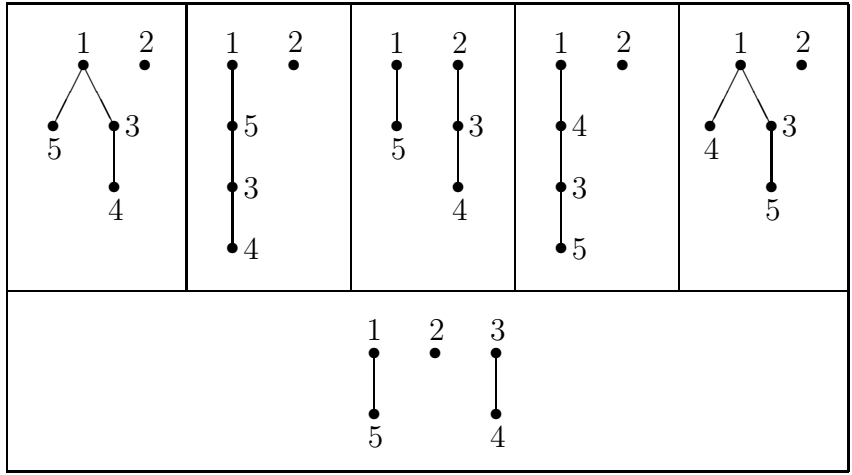

It is clear that
\begin{equation}\label{eq-clear}
|\mathcal{F}_n^{n-1}[n\prec 1]|=1.
\end{equation}
We have the following corollaries.
\begin{cor}{\rm\bf(Cayley \cite{Cayley})}
The number of labeled trees on $n$ vertices is $n^{n-2}$.
\end{cor}
\begin{cor}{\rm\bf(Cayley \cite{Cayley}, Clarke \cite{Clarke})}
The number of rooted trees on $n+1$ vertices with a specific root
and root degree $k$ is ${n-1\choose k-1}n^{n-k}$.
\end{cor}

\noindent{\it Proof.} It follows from \eqref{eq-main} and \eqref{eq-clear} that
\begin{equation*}
|\mathcal{F}_n^k[n\prec 1]|=n^{n-k-1}.
\end{equation*}

Exchanging the labels of the vertices $j$ and $1$ for $1<j\leqslant k<n$,
we establish a bijection between $\mathcal{F}_n^k[n\prec 1]$ and
$\mathcal{F}_n^k[n\prec j]$. Therefore,
\begin{equation*}
|\mathcal{F}_n^k|=kn^{n-k-1},
\end{equation*}
from which one can see that the number of forests with $n$
vertices and $k$ trees is
\begin{equation*}
{n\choose k}|\mathcal{F}_n^k|={n\choose k}kn^{n-k-1}
 ={n-1\choose k-1}n^{n-k}.
\end{equation*}
\qed

\noindent{\it Remark.} It is worth mentioning the third and fourth proofs of Cayley's formula
in Aigner and Ziegler's book \cite[Chapter 30]{AZ}. The third proof in \cite[Chapter 30]{AZ}, essentially
due to Riordan \cite{Riordan} and R\'enyi \cite{Renyi}, is as follows:
Let $T_{n,k}$ denote the number of forests on $[n]$ consisting $k$ trees
where the vertices of $[k]$ appear in different trees. Consider such a forest $F$
and suppose that $1$ is adjacent to $i$ vertices. Removing the vertex $1$, we obtain a forest
of $k-1+i$ trees. As we can reconstruct $F$ by first fixing $i$, then selecting the
$i$ neighbors of $1$, and then the forest $F\setminus 1$, this gives
$$
T_{n,k}=\sum_{i=0}^{n-k}{n-k\choose i}T_{n-1,k-1+i},
$$
from which we can prove Cayley's formula by induction on $n$.

The fourth proof in \cite[Chapter 30]{AZ},
due to Pitman \cite{Pitman}, is as follows:
Let $\mathcal{F}_{n,k}$ denote the set of all forests that consist of $k$
rooted trees on $[n]$. A sequence $F_1,\ldots,F_k$ of forests is called a
refining sequence if $F_i\in\mathcal{F}_{n,i}$ and $F_{i+1}$ is obtained from $F_{i}$ by removing one edge for all $i$.
Fix $F_k\in\mathcal{F}_{n,k}$ and denote by
\begin{itemize}
\item $N(F_k)$ the number of rooted trees containing $F_k$,
\item $N^*(F_k)$ the number of refining sequences ending in $F_k$
\end{itemize}
Then it is easy to see that $N^*(F_k)=N(F_k)(k-1)!$ and $N^*(F_k)=n^{k-1}(k-1)!$
by counting $N^*(F_k)$ in two ways, first by starting a tree and secondly by starting at $F_k$.
Hence,
$$N(F_k)=n^{k-1}\quad\text{for any $F_k\in\mathcal{F}_{n,k}$.}$$
Since $F_n$ contains $n$ isolated nodes, we get $\mathcal{F}_{n,1}=n^{n-1}$, and thus Cayley's formula.

Therefore, both our proof of Cayley's formula and the third proof in \cite[Chapter 30]{AZ} are recursive.
But the latter requires more algebraic computations.
Our proof is very similar to Pitman's proof. The difference is that Pitman
considers all forests with $k$ rooted trees while we only consider forests with roots $1,\dots,k$. Pitman's proof
uses the idea of double counting and our proof is a little more straightforward and combinatorial.

\section{Bipartite trees and tripartite trees}
A complete $k$-partite graph $K_{n_1,\ldots,n_k}$ is a graph
$G=(V,E)$ with the vertex set $V$ partitioned into $k$ parts
$V_1,\ldots, V_k$, where $|V_i|=n_i,\,1\leqslant i\leqslant k$, any two
vertices $x,y\in V$ are adjacent if and only if $x$ and $y$ do not
belong to the same part $V_i$. A spanning tree of
$K_{n_1,\ldots,n_k}$ is called a $k$-partite tree.

\begin{cor}\label{cor-Fiedler}
{\rm\bf(Fiedler and Sedl$\acute{\rm a}\check{\rm c}$ek \cite{Fiedler})}
The number $c(K_{r,s})$ of spanning trees of the complete
bipartite graph $K_{r,s}$ is $r^{s-1}s^{r-1}$.
\end{cor}

\noindent{\it Proof.} Let $V_1=\{1,\ldots,r\}$ and $V_2=\{r+1,\ldots,r+s\}$.
Denote by $\mathcal{F}_{r,s}^k$ the set of forests of $k$ rooted
bipartite trees on $V=V_1\cup V_2$ with roots $1,\ldots,k$.

Note there is no edge connecting vertices $x$ and $y$ if $x$ and $y$
belong to the same set $V_i$. In a similar way to Theorem
\ref{thm-main}, we have
\begin{equation}\label{eq-fb}
 |\mathcal{F}_{r,s}^{k-1}[r+1\prec 1]|=s|\mathcal{F}_{r,s}^k[r+1\prec 1]|,
 \quad 2\leqslant k\leqslant r.
\end{equation}

It is easy to see that
\begin{equation}\label{fb-easy}
|\mathcal{F}_{r,s}^{r}[r+1\prec 1]|=r^{s-1}.
\end{equation}
We complete the proof by \eqref{eq-fb} and \eqref{fb-easy}.\qed

\begin{cor}\label{cor-rs-id}
For $r\geqslant 2,\,s\geqslant 1$, we have
\begin{equation*}
 \sum_{i=1}^{r-1}\sum_{j=1}^{s}{r-2\choose i-1 }{s-1\choose j-1}
 i^{j-1}j^{i-1}(r-i)^{s-j-1}(s-j)^{r-i-1}=r^{s-1}s^{r-2}.
\end{equation*}
Here we set $0^0=1$.
\end{cor}
{\it Proof.} It follows from \eqref{eq-fb} and \eqref{fb-easy}
that
\begin{equation*}
 |\mathcal{F}_{r,s}^{2}[r+1\prec 1]|=r^{s-1}s^{r-2}.
\end{equation*}
On the other hand, both of the trees in a forest $F$ of
$\mathcal{F}_{r,s}^{2}[r+1\prec 1]$ are bipartite trees. Applying
Corollary \ref{cor-Fiedler}, we complete the proof. \qed

\begin{cor}{\rm\bf(Austin \cite{Austin})}
 The number $c(K_{r,s,t})$ of spanning trees of the complete
 tripartite graph $K_{r,s,t}$ is $(r+s+t)(r+s)^{t-1}(s+t)^{r-1}(t+r)^{s-1}$.
\end{cor}

\noindent{\it Proof.} Let $V_1=\{1,\ldots,r\}$, $V_2=\{r+1,\ldots,r+s\}$, and
$V_3=\{r+s+1,\ldots,r+s+t\}$. Denote by $\mathcal{F}_{r,s,t}^k$
(respectively, $\check{\mathcal{F}}_{r,s,t}^k$) the set of forests of $k$
rooted tripartite trees on $V=V_1\cup V_2\cup V_3$ with roots
$1,\ldots,k$ (respectively, $2,\ldots,k+1$).

In exactly the same way as in the proof of Corollary \ref{cor-Fiedler}, we
have
\begin{equation*}
 |\mathcal{F}_{r,s,t}^{k-1}[r+1\prec 1]|=(s+t)|\mathcal{F}_{r,s,t}^k[r+1\prec 1]|,
 \quad 2\leqslant k\leqslant r.
\end{equation*}
It follows that
\begin{equation}\label{rst-r}
c(K_{r,s,t})=
 |\mathcal{F}_{r,s,t}^{1}[r+1\prec 1]|=(s+t)^{r-1}|\mathcal{F}_{r,s,t}^r[r+1\prec
  1]|.
\end{equation}
We have the following important fact:
\begin{equation}\label{rst-r+1}
 |\mathcal{F}_{r,s,t}^r[r+1\prec 1]|=|\check{\mathcal{F}}_{r,s,t}^r[1\prec r+1]|,
\end{equation}
because $1$ and $r+1$ are in the same tree in each forest, what is
different is that we choose $1$ as the root of the tree on the
left-hand side, while we choose $r+1$ on the right-hand side. Again,
similarly to the proof of Corollary \ref{cor-Fiedler}, we get
\begin{equation}\label{rst-s}
 |\check{\mathcal{F}}_{r,s,t}^{k-1}[1\prec r+1]|=
 (r+t)|\check{\mathcal{F}}_{r,s,t}^k[1\prec r+1]|,
 \quad r+1\leqslant k\leqslant r+s-1.
\end{equation}

On the other hand, we have
\begin{equation}\label{rst-last}
 |\check{\mathcal{F}}_{r,s,t}^{r+s-1}[1\prec r+1]|=
(r+s+t)(r+s)^{t-1}.
\end{equation}
In fact, there are $(r+s)^t$ forests such that $1$ is a child of
$r+1$, and $t(r+s)^{t-1}$ forests such that the father of $1$ is a
child of $r+1$, which comes from $V_3$.

The proof then follows from combining \eqref{rst-r}--\eqref{rst-last}.\qed

\noindent{\it Remark.} A bijective proof of Corollary \ref{cor-Fiedler} was found by Stanley \cite[pp.~125-126]{Stanley}.
Austin \cite{Austin} proved that the number of spanning trees of the complete
$k$-partite graph $K_{n_1,\ldots,n_k}$ is $n^{k-2}(n-n_1)^{n_1-1}\cdots(n-n_k)^{n_k-1}$, where $n=n_1+\cdots+n_p$.
E\v{g}ecio\v{g}lu and Remmel\cite{ER86} found a bijective proof of the enumeration of bipartite trees
and tripartite trees and later they \cite{ER94} found a bijection for the $k$-partite trees.
The reader is encouraged to find a simple proof of Austin's formula by modifying our recursive algorithm for
trees and forests.

\section{Plane trees}

A labeled plane tree is a rooted labeled tree for which the
children of any vertex are linearly ordered. Denote by
$\mathcal{P}_{n}^k$ the set of forests of $k$ labeled plane trees
on $n$ vertices with roots $1,\ldots,k$.

\begin{cor}\label{cor-plane}
For $2\leqslant k\leqslant n-1$, the following holds
\begin{equation}\label{eq-plane}
 |\mathcal{P}_n^{k-1}[n\prec 1]|=(2n-k)|\mathcal{P}_n^{k} [n\prec 1]|.
\end{equation}
\end{cor}

\noindent{\it Proof.} The first step is the same as Theorem \ref{thm-main}.
For the reversal step, we suppose that $F'\in\mathcal{P}_n^k[n\prec 1]$.
For the first case, there are $\deg (x)+1$ ways for us to attach
$F'_k$ to the vertex $x$ of the other trees in $F$ as a subtree. The second case is the same. The proof then follows from
the fact
\begin{equation*}
\sum_{x\in V(F')}(\deg_{F'}(x)+1)=n-k+n=2n-k.
\end{equation*}

\begin{cor}
The number of labeled plane trees on $n+1$ vertices is $(2n)!/n!$,
i.e., the number of unlabeled plane trees on $n+1$ vertices is the
well-known Catalan number
\begin{equation*}
n_n=\frac{1}{n+1}{2n\choose n}.
\end{equation*}
\end{cor}

\noindent{\it Proof.} It follows from \eqref{eq-plane} and
$|\mathcal{P}_{n+1}^n[n+1\prec 1]|=1$ that the number of labeled
plane trees on $[n+1]$ with root $1$ is $(2n)!/(n+1)!$, but the
total number of plane trees is $n+1$ times that.\qed

We say that a plane tree is {\it non-leaf vertex labeled} if all the leaves are not labeled, while
all of the other vertices are labeled.
Denote by $\mathcal{P}_{n,p}^r$ the set of forests of $k$ non-leaf vertex labeled plane trees
on $n$ vertices with $p$ leaves and with roots $1,\ldots,r$.

\begin{thm}\label{thm-p-leaf}
For $2\leqslant r<n-p$, the following holds
\begin{equation}\label{eq-p-leaf}
 |\mathcal{P}_{n,p}^{r-1}[n-p\prec 1]|=(p+1)|\mathcal{P}_{n+1,p+1}^{r} [n-p\prec 1]|.
\end{equation}
\end{thm}

\noindent{\it Proof.} The proof is analogous to that of Theorem
\ref{thm-main}. But this time when a subtree is removed, its
root is left and is considered as an unlabeled vertex. Thus, the number of the
leaves (unlabeled vertices) of the forest is increased by one. The
reverse is clear and straightforward. \qed

\begin{figure}[h]
\begin{center}
\setlength{\unitlength}{0.08in}
\begin{tabular}{|c|}\hline
\begin{picture}(14,18)
\put(7,14){\line(-1,-2){6}}
\put(7,14){\line(1,-2){2}}
\put(5,10){\line(1,-2){2}}
\put(3,6){\line(1,-2){2}}
\put(3,6){\line(0,-1){4}}
\put(12,14){\line(0,-1){8}}
\put(7,14){\circle*{0.5}}    \put(6.5,14.7){\makebox(2,1)[l]{1}}
\put(5,10){\circle*{0.5}}    \put(3.6,10.2){\makebox(2,1)[l]{5}}
\put(3,6){\circle*{0.5}}     \put(1.6,6.2){\makebox(2,1)[l]{3}}
\put(1,2){\circle*{0.5}}
\put(3,2){\circle*{0.5}}
\put(5,2){\circle*{0.5}}
\put(3,6){\circle*{0.5}}
\put(7,6){\circle*{0.5}}
\put(9,10){\circle*{0.5}}
\put(12,14){\circle*{0.5}}    \put(11.6,14.7){\makebox(2,1)[l]{2}}
\put(12,10){\circle*{0.5}}    \put(12.6,9.5){\makebox(2,1)[l]{4}}
\put(12,6){\circle*{0.5}}
\put(7,14){\line(0,-1){4}}\put(7,10){\circle*{0.5}}
\end{picture}
\vline
\begin{picture}(14,18)
\put(7,14){\line(-1,-2){4}}
\put(7,14){\line(1,-2){2}}
\put(5,10){\line(1,-2){4}}
\put(7,6){\line(-1,-2){2}}
\put(7,6){\line(0,-1){4}}
\put(12,14){\line(0,-1){8}}
\put(7,14){\circle*{0.5}}    \put(6.5,14.7){\makebox(2,1)[l]{1}}
\put(9,10){\circle*{0.5}}
\put(5,10){\circle*{0.5}}    \put(3.6,10.2){\makebox(2,1)[l]{5}}
\put(3,6){\circle*{0.5}}
\put(7,6){\circle*{0.5}}     \put(7.4,6.2){\makebox(2,1)[l]{3}}
\put(5,2){\circle*{0.5}}
\put(7,2){\circle*{0.5}}
\put(9,2){\circle*{0.5}}
\put(12,14){\circle*{0.5}}    \put(11.6,14.7){\makebox(2,1)[l]{2}}
\put(12,10){\circle*{0.5}}    \put(12.6,9.5){\makebox(2,1)[l]{4}}
\put(12,6){\circle*{0.5}}
\put(7,14){\line(0,-1){4}}\put(7,10){\circle*{0.5}}
\end{picture}
\vline
\begin{picture}(14,18)
\put(6,14){\line(-1,-2){4}}
\put(6,14){\line(3,-4){3}}
\put(6,10){\line(1,-2){2}}
\put(6,10){\line(1,-1){4}}
\put(4,10){\line(0,-1){4}}
\put(12,14){\line(0,-1){8}}
\put(6,14){\circle*{0.5}}    \put(5.5,14.7){\makebox(2,1)[l]{1}}
\put(4,10){\circle*{0.5}}    \put(2.6,10.2){\makebox(2,1)[l]{5}}
\put(2,6){\circle*{0.5}}
\put(4,6){\circle*{0.5}}
\put(9,10){\circle*{0.5}}
\put(6,6){\circle*{0.5}}
\put(8,6){\circle*{0.5}}
\put(10,6){\circle*{0.5}}
\put(12,14){\circle*{0.5}}   \put(11.6,14.7){\makebox(2,1)[l]{2}}
\put(12,10){\circle*{0.5}}   \put(12.6,9.5){\makebox(2,1)[l]{4}}
\put(12,6){\circle*{0.5}}
\put(6,14){\line(0,-1){8}}    \put(6,10){\circle*{0.5}}
                              \put(6.4,10.2){\makebox(2,1)[l]{3}}
\end{picture}
   \vline
\begin{picture}(14,18)
\put(6,14){\line(-1,-2){4}}
\put(6,14){\line(1,-2){4}}
\put(4,10){\line(0,-1){4}}
\put(8,10){\line(-1,-2){2}}
\put(8,10){\line(0,-1){4}}
\put(12,14){\line(0,-1){8}}
\put(6,14){\circle*{0.5}}    \put(5.5,14.7){\makebox(2,1)[l]{1}}
\put(4,10){\circle*{0.5}}    \put(2.6,10.2){\makebox(2,1)[l]{5}}
\put(2,6){\circle*{0.5}}
\put(4,6){\circle*{0.5}}
\put(8,10){\circle*{0.5}}     \put(8.4,10.2){\makebox(2,1)[l]{3}}
\put(6,6){\circle*{0.5}}
\put(8,6){\circle*{0.5}}
\put(10,6){\circle*{0.5}}
\put(12,14){\circle*{0.5}}   \put(11.6,14.7){\makebox(2,1)[l]{2}}
\put(12,10){\circle*{0.5}}   \put(12.6,9.5){\makebox(2,1)[l]{4}}
\put(12,6){\circle*{0.5}}
\put(6,14){\line(0,-1){4}}\put(6,10){\circle*{0.5}}
\end{picture}
\\ \hline
\begin{picture}(14,18)
\put(6,14){\line(-1,-2){4}}
\put(6,14){\line(1,-2){2}}
\put(4,10){\line(1,-2){2}}
\put(11,14){\line(0,-1){12}}
\put(11,6){\line(-1,-2){2}}
\put(11,6){\line(1,-2){2}}
\put(6,14){\circle*{0.5}}    \put(5.5,14.7){\makebox(2,1)[l]{1}}
\put(4,10){\circle*{0.5}}    \put(2.6,10.2){\makebox(2,1)[l]{5}}
\put(2,6){\circle*{0.5}}
\put(6,6){\circle*{0.5}}
\put(8,10){\circle*{0.5}}
\put(11,14){\circle*{0.5}}   \put(10.6,14.7){\makebox(2,1)[l]{2}}
\put(11,10){\circle*{0.5}}   \put(11.6,9.5){\makebox(2,1)[l]{4}}
\put(11,6){\circle*{0.5}}    \put(11.6,5.5){\makebox(2,1)[l]{3}}
\put(9,2){\circle*{0.5}}
\put(11,2){\circle*{0.5}}
\put(13,2){\circle*{0.5}}
\put(6,14){\line(0,-1){4}}\put(6,10){\circle*{0.5}}
\end{picture}
\vline
\begin{picture}(14,18)
\put(8,14){\line(-1,-2){6}}
\put(8,14){\line(0,-1){4}}
\put(8,14){\line(1,-2){2}}
\put(6,10){\line(1,-2){2}}
\put(4,6){\line(1,-2){2}}
\put(12,14){\line(0,-1){8}}
\put(8,14){\circle*{0.5}}    \put(7.5,14.7){\makebox(2,1)[l]{1}}
\put(8,10){\circle*{0.5}}
\put(6,10){\circle*{0.5}}    \put(4.6,10.2){\makebox(2,1)[l]{3}}
\put(4,6){\circle*{0.5}}    \put(2.6,6.2){\makebox(2,1)[l]{5}}
\put(2,2){\circle*{0.5}}
\put(6,2){\circle*{0.5}}
\put(8,6){\circle*{0.5}}
\put(10,10){\circle*{0.5}}
\put(12,14){\circle*{0.5}}    \put(11.6,14.7){\makebox(2,1)[l]{2}}
\put(12,10){\circle*{0.5}}     \put(12.6,9.5){\makebox(2,1)[l]{4}}
\put(12,6){\circle*{0.5}}
\put(6,10){\line(0,-1){4}}\put(6,6){\circle*{0.5}}
\end{picture}
\vline
\begin{picture}(14,18)
\put(6,14){\line(-1,-2){2}}
\put(6,14){\line(0,-1){4}}
\put(6,10){\line(-1,-2){4}}
\put(6,14){\line(3,-4){3}}
\put(6,10){\line(1,-2){2}}
\put(4,6){\line(1,-2){2}}
\put(12,14){\line(0,-1){8}}
\put(6,14){\circle*{0.5}}    \put(5.5,14.7){\makebox(2,1)[l]{1}}
\put(4,10){\circle*{0.5}}
\put(4,6){\circle*{0.5}}     \put(2.6,6.2){\makebox(2,1)[l]{5}}
\put(2,2){\circle*{0.5}}
\put(6,2){\circle*{0.5}}
\put(8,6){\circle*{0.5}}
\put(6,10){\circle*{0.5}}    \put(6.4,10.2){\makebox(2,1)[l]{3}}
\put(9,10){\circle*{0.5}}
\put(12,14){\circle*{0.5}}    \put(11.6,14.7){\makebox(2,1)[l]{2}}
\put(12,10){\circle*{0.5}}    \put(12.6,9.5){\makebox(2,1)[l]{4}}
\put(12,6){\circle*{0.5}}
\put(6,10){\line(0,-1){4}}\put(6,6){\circle*{0.5}}
\end{picture}
   \vline
\begin{picture}(14,18)
\put(5,14){\line(-1,-2){2}}
\put(5,14){\line(0,-1){4}}
\put(5,14){\line(1,-2){4}}
\put(7,10){\line(-1,-2){4}}
\put(5,6){\line(1,-2){2}}
\put(12,14){\line(0,-1){8}}
\put(5,14){\circle*{0.5}}    \put(4.5,15){\makebox(2,1)[l]{1}}
\put(3,10){\circle*{0.5}}
\put(5,10){\circle*{0.5}}
\put(7,10){\circle*{0.5}}   \put(7.4,10.2){\makebox(2,1)[l]{3}}
\put(5,6){\circle*{0.5}}    \put(3.6,6.2){\makebox(2,1)[l]{5}}
\put(3,2){\circle*{0.5}}
\put(7,2){\circle*{0.5}}
\put(9,6){\circle*{0.5}}
\put(12,14){\circle*{0.5}}   \put(11.6,14.7){\makebox(2,1)[l]{2}}
\put(12,10){\circle*{0.5}}   \put(12.6,9.5){\makebox(2,1)[l]{4}}
\put(12,6){\circle*{0.5}}
\put(7,10){\line(0,-1){4}}\put(7,6){\circle*{0.5}}
\end{picture}
\\ \hline
\begin{picture}(27,14)
\put(8,10){\line(-1,-2){4}}
\put(20,10){\line(0,-1){4}}
\put(8,10){\line(1,-2){2}}
\put(6,6){\line(1,-2){2}}
\put(14,10){\line(0,-1){8}}
\put(20,10){\line(-1,-2){2}}
\put(20,10){\line(1,-2){2}}
\put(8,10){\circle*{0.5}}    \put(7.5,10.7){\makebox(2,1)[l]{1}}
\put(6,6){\circle*{0.5}}     \put(4.6,6.2){\makebox(2,1)[l]{5}}
\put(4,2){\circle*{0.5}}
\put(8,2){\circle*{0.5}}
\put(10,6){\circle*{0.5}}
\put(14,10){\circle*{0.5}}    \put(13.6,10.7){\makebox(2,1)[l]{2}}
\put(14,6){\circle*{0.5}}     \put(14.6,5.5){\makebox(2,1)[l]{4}}
\put(14,2){\circle*{0.5}}
\put(20,10){\circle*{0.5}}    \put(19.6,10.7){\makebox(2,1)[l]{3}}
\put(18,6){\circle*{0.5}}
\put(20,6){\circle*{0.5}}
\put(22,6){\circle*{0.5}}
\put(8,10){\line(0,-1){4}}\put(8,6){\circle*{0.5}}
\end{picture}
\\ \hline
\end{tabular}
\caption{An example for Theorem \ref{thm-p-leaf}.}
\label{fig-plane}
\end{center}
\end{figure}
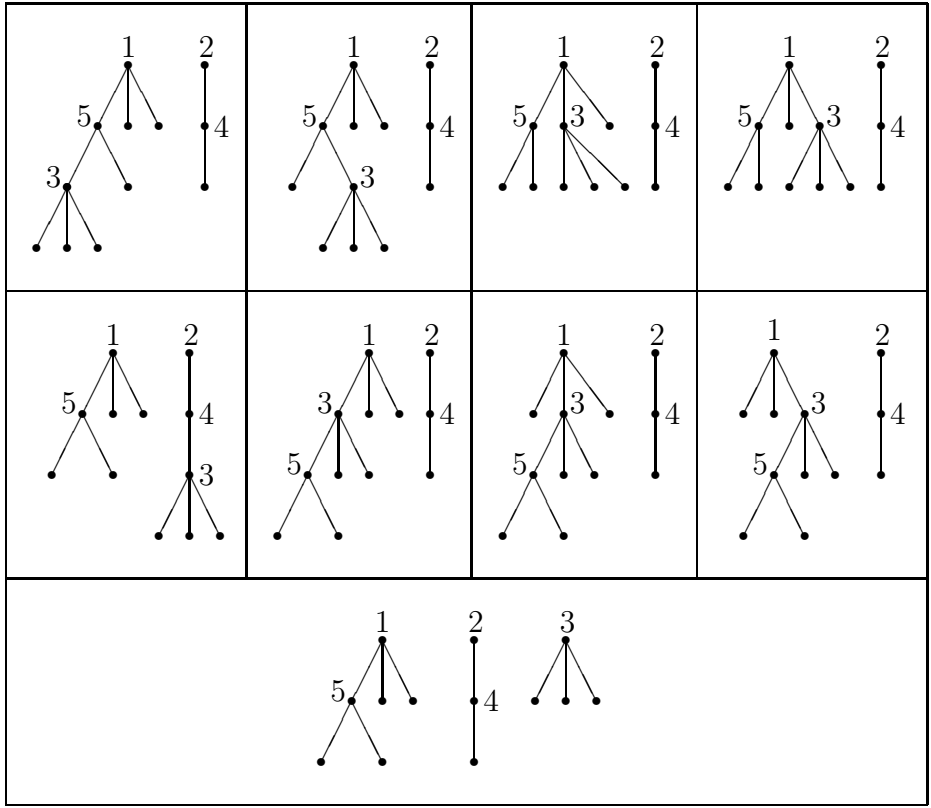

For $n=12,\,p=7,\,r=3$, an example is given in Figure \ref{fig-plane}.

A plane tree is called a $k$-ary plane tree if every non-leaf vertex has degree $k$.
\begin{cor}\label{cor-k-ary}
The number of forests of $r$ non-leaf vertex labeled $k$-ary plane trees on $kn+r$ vertices
with roots $1,\ldots,r$ equals
\[
\frac{r}{n}{kn\choose n-r}(n-r)!.
\]
\end{cor}

\noindent{\it Proof.} The recurrence \eqref{eq-p-leaf} also holds for $k$-ary plane trees.
But in this case, $(n,p)$ is replaced by $(kn+1, (k-1)n+1)$. Hence, we have
\[
|\mathcal{P}_{kn+1,(k-1)n+r}^{r}[n\prec 1]|=
((k-1)n+r+1)((k-1)n+r+2)\cdots(kn-1)\cdot k.
\]
Here the factor $k$ occurs because there are $k$ possibilities for $n$ to be a child of $1$.
The proof then follows from the fact that
$|\mathcal{P}_{kn+1,(k-1)n+1}^{r}|=r|\mathcal{P}_{kn+1,(k-1)n+1}^{r}[n\prec 1]|$. \qed

\begin{cor}\label{cor-ary-tree}
The number of unlabeled $k$-ary plane trees on $kn+1$ vertices equals
\[
\frac{1}{n}{kn\choose n-1}=\frac{1}{kn+1}{kn+1\choose n}.
\]
\end{cor}

\noindent{\it Proof.} For any $k$-ary plane tree on $kn+1$ vertices, it has $n$ non-leaf vertices.
We have $(n-1)!$ ways to label its non-leaf vertices with the root labeled $1$.
Therefore, dividing the number of non-leaf labeled plane trees on $kn+1$ vertices with
root $1$ by $(n-1)!$, we obtain the desired number for unlabeled $k$-ary plane trees. \qed

\begin{cor}\label{cor-ary-id}
For $p,\,q\geqslant 1,\,n\geqslant p+q$, we have
\begin{equation*}
\sum_{i=p}^{n-q}\frac{pq}{i(n-i)}{ki\choose i-p}{k(n-i)\choose n-i-q}
=\frac{p+q}{n}{kn\choose n-(p+q)}.
\end{equation*}
\end{cor}

\noindent{\it Proof.} Let $r=p+q$. The forest of $p+q$ non-leaf labeled trees may be treated as
the union of two small forests, namely, the forest of $p$ trees with roots $1,\ldots, p$,
and the forest of $q$ trees with roots $p+1,\ldots,p+q$. It follows from
Corollary~\ref{cor-k-ary} that
\begin{align*}
&\hskip -2mm
\sum_{i=p}^{n-q}{n-(p+q)\choose i-p}\frac{p}{i}{ki\choose i-p}(i-p)!
\frac{q}{n-i}{k(n-i)\choose n-i-q}(n-i-q)! \\[5pt]
&=\frac{p+q}{n}{kn\choose n-(p+q)}(n-(p+q))!.
\end{align*}
Dividing both sides of the identity by $(n-(p+q))!$, we complete the proof.
\qed

\begin{cor}{\rm\bf(Narayana \cite{Narayana})}\label{cor-Nara}
The number of unlabeled plane trees on $n+1$ vertices with $p$
leaves equals
\[
\frac{1}{n+1}{n+1\choose p}{n-1\choose n-p}.
\]
\end{cor}

\noindent{\it Proof.} It follows from \eqref{eq-p-leaf} that
\begin{align}\label{eq-factor}
\nonumber
 |\mathcal{P}_{n+1,p}^{1}[n+1-p\prec 1]|
 &=(p+1)(p+2)\cdots(n-1)|\mathcal{P}_{2n-p,n-1}^{n-p} [n+1-p\prec 1]|\nonumber\\[5pt]
 &=\frac{(n-1)!}{p!}|\mathcal{P}_{2n-p,n-1}^{n-p} [n+1-p\prec 1]|.
\end{align}

For any $F\in\mathcal{P}_{2n-p,n-1}^{n-p} [n+1-p\prec 1]$, one sees that $n+1-p$ is a
child of $1$ and so the number of forests in $\mathcal{P}_{2n-p,n-1}^{n-p} [n+1-p\prec 1]$ with the vertex $i$ having degree $x_i$
is equal to $x_1$ (there are $x_1$ possibilities for $n+1-p$ to be a child of $1$). It follows that
\begin{equation}\label{eq-sum}
|\mathcal{P}_{2n-p,n-1}^{n-p} [n+1-p\prec 1]|=\sum_{x_1+\cdots+x_{n-p+1}=n} x_1.
\end{equation}

It is well known that the number of solutions in positive integers to the equation
\begin{equation*}
x_1+\cdots+x_{n-p+1}=n,
\end{equation*}
is ${n-1\choose n-p}$. So the right-hand side of \eqref{eq-sum} is equal to
\[
\frac{n}{n-p+1}{n-1\choose n-p}.
\]
Together
with \eqref{eq-factor}, we obtain
\begin{equation*}
 |\mathcal{P}_{n+1,p}^{1}|
 =\frac{(n-p)!}{n+1}{n+1\choose p}{n-1\choose n-p}.
\end{equation*}
Dividing this number by $(n-p)!$, we get the desired formula for
unlabeled plane trees.\qed

\begin{cor}\label{cor-p-deg}
For any nonnegative integer sequence $(d_1,d_2,\ldots,d_n)$ satisfying
$d_1+d_2+\cdots+d_n=n-1$, there are $(n-1)!$ plane trees on $[n]$
with the degree of $i$ being $d_i$.
\end{cor}

\noindent{\it Proof.} Without loss of generality, suppose that $d_i>0$ if and only if
$1\leqslant i\leqslant n-p$. Just like \eqref{eq-factor}, we see that the number of non-leaf
vertex labeled plane trees with root $i$ and with the degree of $j$ being $d_j$ for
$1\leqslant j\leqslant n-p$ is $d_i(n-2)!/p!$, so the total number is
\[
\sum_{d_i>0}d_i\frac{(n-2)!}{p!}=\frac{(n-1)!}{p!}.
\]
But there are $p!$ ways to label the leaves of each plane tree with
$\{n-p+1,\ldots,n\}$. This completes the proof.\qed

\begin{cor}\label{cor-tree-deg}
For any nonnegative integer sequence $(d_1,d_2,\ldots,d_n)$ satisfying
$d_1+d_2+\cdots+d_n=n-1$, the number of rooted trees on
$[n]$ with the degree of $i$ being $d_i$ is
\[
{n-1\choose d_1,d_2,\ldots,d_n}.
\]
\end{cor}

\noindent{\it Proof.} For any rooted tree on $[n]$ with the degree of $i$ being $d_i$, there are
$d_1!d_2!\ldots d_n!$ plane trees with the same degree sequence.
The proof then follows from Corollary \ref{cor-p-deg}.\qed

\begin{cor}{\rm\bf(Erd$\acute{\rm e}$lyi-Etherington \cite{EE})}\label{cor-EE}
Let $1^{n_1}2^{n_2}\cdots m^{n_m}$ be a partition of $n-1$,
i.e., $n_1+2n_2+\cdots+mn_m=n-1$, and $n_0=n-(n_1+n_2+\cdots+n_m)$. Then the number
of unlabeled plane trees having $n_i$ vertices with degree $i$ is equal to
\[
\frac{1}{n}{n\choose{n_0,n_1,\ldots,n_m}}.
\]
\end{cor}

\noindent{\it Proof.} It follows from the proof of Corollary \ref{cor-Nara} (replace $n$ by $n-1$)
that the number of such non-leaf vertex labeled plane trees is equal to the product of
$(n-2)!/n_0!$ and  $\sum x_1$, where $x_1,x_2,\ldots,x_{n-n_0}$ ranges over all the
permutations of the multiset $\{1^{n_1},2^{n_2},\ldots,m^{n_m}\}$. It is not hard to see
that
\begin{equation*}
\sum x_1=\sum_{j=1}^{m}\frac{(n_1+\cdots+n_m-1)!}{n_1!\cdots n_m!/{n_j}}j=
\frac{(n_1+\cdots+n_m-1)!}{n_1!\cdots n_m!}(n-1).
\end{equation*}
That is, the multiplication equals
\[
\frac{(n_1+\cdots+n_m-1)!}{n_0!n_1!\cdots n_m!}(n-1)!.
\]
Dividing it by $(n-1-n_0)!$, we obtain the desired number for unlabeled plane trees.
\qed

\section{$k$-Edge colored trees}

A {\it $k$-edge colored tree} is a tree whose edges
are colored from a set of $k$ colors such that any two edges with a common
vertex have different colors. For convenience, the set of $k$ colors is denoted by
$\{{\bf 1},{\bf 2},\ldots,{\bf k}\}$.

A {\it special $k$-edge colored tree} is a rooted $k$-edge colored tree such that
all the edges out of the root are colored with the first $k-1$ colors, i.e.,
$\{{\bf 1},{\bf 2},\ldots,{\bf k-1}\}$. Note that for such a tree, the degree of
the root is less than or equal to $k-1$. Denote by $\mathcal{E}_{n,k}^{r}$ the set of forests of
$r$ special $k$-edge trees with roots $1,\ldots,r$.

\begin{thm}\label{thm-color}
For $2\leqslant r\leqslant n-1$, we have the following recurrence relations:
\begin{equation}
\label{eq-special}
|\mathcal{E}_{n,k}^{r-1}[n\prec 1]|=(kn-2n+r)|\mathcal{E}_{n,k}^{r}[n\prec 1]|.
\end{equation}
\end{thm}

\noindent{\it Proof.} For $F\in\mathcal{E}_{n,k}^{r-1}[n\prec 1]$, assume that the edge with child
vertex $r$ has the color ${\bf x}$, then none of the edges with father vertex $r$
could have color ${\bf x}$.
\begin{itemize}
\item First, remove the subtree $F_r$ from $F$ and add it to be a new tree in the forest.
\item
Second, if one edge with father vertex $r$ is colored by ${\bf k}$,
we recolor it by ${\bf x}$.
\item
Third, if $n$ is not a descendant of $1$ in the new forest, $n$ must be a descendant
of $r$, exchange the labels of vertices $1$ and $r$.
\end{itemize}
Thus we obtain a forest $F'\in\mathcal{E}_{n,k}^{r}[n\prec 1]$.

Conversely, for a forest $F'\in\mathcal{E}_{n,k}^{r}[n\prec 1]$, we can use the
following steps to recover $F$:

\begin{itemize}
\item First, attach $F'_r$ back to any vertex of the other trees in $F$
as a subtree, or attach $F'_1$ to any vertex of $F'_r$ as a subtree and
exchange the labels of the vertices $1$ and $r$. Note that there is an uncolored edge.
i.e., the edge $e$ with child vertex $r$. Denote by $v$ the father vertex of $e$,
and $C(v)$ the set of colors of the other edges incident with $v$.

\item Second, if $v$ is one of the roots in the forest, i.e., $1\leqslant v\leqslant r$, we can color
 $e$ with any color in $\{{\bf 1},{\bf 2},\ldots,{\bf k-1}\}\setminus C(v)$,
 whenever it is not empty.  Otherwise, we color $e$ with any color in
 $\{{\bf 1},{\bf 2},\ldots,{\bf k}\}\setminus C(v)$, whenever it is not empty.

\item Third, if the color of $e$, say ${\bf x}$, is the same as some edge $f$
incident with vertex $r$, then recolor $f$ by ${\bf k}$.
\end{itemize}

Note that if $1\leqslant v\leqslant r$, we have $C(v)\subseteq\{{\bf 1},{\bf 2},\ldots,{\bf k-1}\}$,
and $|C(v)|=\deg_{F'}(v)$. Otherwise, we have $|C(v)|=\deg_{F'}(v)+1$. Therefore,
the number of ways to recover $F'$ is given by
\begin{align*}
\nonumber
&\hskip -2mm
\sum_{v=1}^{r}(k-1-\deg_{F'}(v))+\sum_{v=r+1}^{n}(k-(1+\deg_{F'}(v)))\\[5pt]
&=(k-1)n-\sum_{v=1}^{n}\deg_{F'}(v)=kn-2n+r.
\end{align*}
This completes the proof of \eqref{eq-special}.\qed

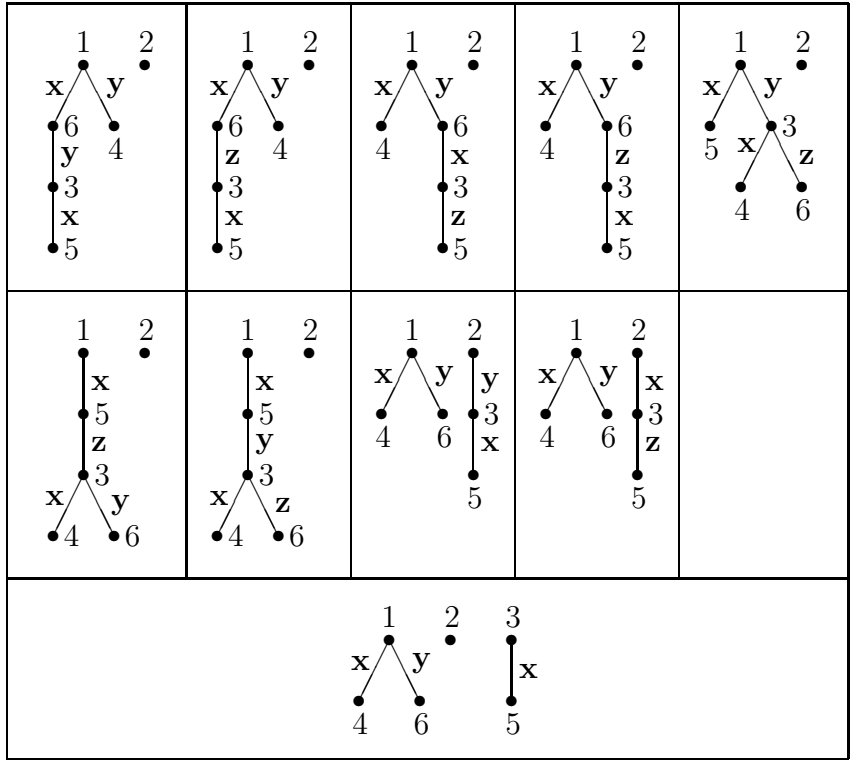
\begin{figure}[h]
\begin{center}
\setlength{\unitlength}{0.08in}
\begin{tabular}{|c|}\hline
\begin{picture}(10,18)
\put(4,14){\line(-1,-2){2}}
\put(4,14){\line(1,-2){2}}
\put(2,10){\line(0,-1){8}}
\put(4,14){\circle*{0.7}}    \put(3.5,15){\makebox(2,1)[l]{1}}
\put(2,10){\circle*{0.7}}    \put(2.7,9.5){\makebox(2,1)[l]{6}}
\put(2,6){\circle*{0.7}}     \put(2.7,5.5){\makebox(2,1)[l]{3}}
\put(2,2){\circle*{0.7}}     \put(2.7,1,5){\makebox(2,1)[l]{5}}
\put(6,10){\circle*{0.7}}    \put(5.6,8){\makebox(2,1)[l]{4}}
\put(8,14){\circle*{0.7}}    \put(7.6,15){\makebox(2,1)[l]{2}}

\put(1.5,12){\makebox(2,1)[l]{\bf x}}
\put(5.5,12){\makebox(2,1)[l]{\bf y}}
\put(2.5,7.5){\makebox(2,1)[l]{\bf y}}
\put(2.5,3.5){\makebox(2,1)[l]{\bf x}}
\end{picture}
   \vline
\begin{picture}(10,18)
\put(4,14){\line(-1,-2){2}}
\put(4,14){\line(1,-2){2}}
\put(2,10){\line(0,-1){8}}
\put(4,14){\circle*{0.7}}    \put(3.5,15){\makebox(2,1)[l]{1}}
\put(2,10){\circle*{0.7}}    \put(2.7,9.5){\makebox(2,1)[l]{6}}
\put(2,6){\circle*{0.7}}     \put(2.7,5.5){\makebox(2,1)[l]{3}}
\put(2,2){\circle*{0.7}}     \put(2.7,1,5){\makebox(2,1)[l]{5}}
\put(6,10){\circle*{0.7}}    \put(5.6,8){\makebox(2,1)[l]{4}}
\put(8,14){\circle*{0.7}}    \put(7.6,15){\makebox(2,1)[l]{2}}

\put(1.5,12){\makebox(2,1)[l]{\bf x}}
\put(5.5,12){\makebox(2,1)[l]{\bf y}}
\put(2.5,7.5){\makebox(2,1)[l]{\bf z}}
\put(2.5,3.5){\makebox(2,1)[l]{\bf x}}
\end{picture}
   \vline

\begin{picture}(10,18)
\put(4,14){\line(-1,-2){2}}
\put(4,14){\line(1,-2){2}}
\put(6,10){\line(0,-1){8}}
\put(4,14){\circle*{0.7}}    \put(3.5,15){\makebox(2,1)[l]{1}}
\put(2,10){\circle*{0.7}}    \put(1.6,8){\makebox(2,1)[l]{4}}
\put(6,10){\circle*{0.7}}    \put(6.7,9.5){\makebox(2,1)[l]{6}}
\put(6,6){\circle*{0.7}}     \put(6.7,5.5){\makebox(2,1)[l]{3}}
\put(6,2){\circle*{0.7}}     \put(6.7,1,5){\makebox(2,1)[l]{5}}
\put(8,14){\circle*{0.7}}    \put(7.6,15){\makebox(2,1)[l]{2}}

\put(1.5,12){\makebox(2,1)[l]{\bf x}}
\put(5.5,12){\makebox(2,1)[l]{\bf y}}
\put(6.5,7.5){\makebox(2,1)[l]{\bf x}}
\put(6.5,3.5){\makebox(2,1)[l]{\bf z}}
\end{picture}
   \vline

\begin{picture}(10,18)
\put(4,14){\line(-1,-2){2}}
\put(4,14){\line(1,-2){2}}
\put(6,10){\line(0,-1){8}}
\put(4,14){\circle*{0.7}}    \put(3.5,15){\makebox(2,1)[l]{1}}
\put(2,10){\circle*{0.7}}    \put(1.6,8){\makebox(2,1)[l]{4}}
\put(6,10){\circle*{0.7}}    \put(6.7,9.5){\makebox(2,1)[l]{6}}
\put(6,6){\circle*{0.7}}     \put(6.7,5.5){\makebox(2,1)[l]{3}}
\put(6,2){\circle*{0.7}}     \put(6.7,1,5){\makebox(2,1)[l]{5}}
\put(8,14){\circle*{0.7}}    \put(7.6,15){\makebox(2,1)[l]{2}}

\put(1.5,12){\makebox(2,1)[l]{\bf x}}
\put(5.5,12){\makebox(2,1)[l]{\bf y}}
\put(6.5,7.5){\makebox(2,1)[l]{\bf z}}
\put(6.5,3.5){\makebox(2,1)[l]{\bf x}}
\end{picture}
   \vline
\begin{picture}(10,18)
\put(4,14){\line(-1,-2){2}}
\put(4,14){\line(1,-2){4}}
\put(6,10){\line(-1,-2){2}}
\put(4,14){\circle*{0.7}}    \put(3.5,15){\makebox(2,1)[l]{1}}
\put(2,10){\circle*{0.7}}    \put(1.6,8){\makebox(2,1)[l]{5}}
\put(6,10){\circle*{0.7}}    \put(6.7,9.5){\makebox(2,1)[l]{3}}
\put(4,6){\circle*{0.7}}     \put(3.6,4){\makebox(2,1)[l]{4}}
\put(8,6){\circle*{0.7}}     \put(7.6,4){\makebox(2,1)[l]{6}}
\put(8,14){\circle*{0.7}}    \put(7.6,15){\makebox(2,1)[l]{2}}
\put(1.5,12){\makebox(2,1)[l]{\bf x}}
\put(5.5,12){\makebox(2,1)[l]{\bf y}}
\put(3.8,8.3){\makebox(2,1)[l]{\bf x}}
\put(7.8,7.5){\makebox(2,1)[l]{\bf z}}
\end{picture}
\\ \hline
\begin{picture}(10,18)
\put(4,14){\line(0,-1){8}}
\put(4,6){\line(-1,-2){2}}
\put(4,6){\line(1,-2){2}}
\put(4,14){\circle*{0.7}}    \put(3.5,15){\makebox(2,1)[l]{1}}
\put(4,10){\circle*{0.7}}    \put(4.7,9.5){\makebox(2,1)[l]{5}}
\put(4,6){\circle*{0.7}}     \put(4.7,5.5){\makebox(2,1)[l]{3}}
\put(2,2){\circle*{0.7}}     \put(2.7,1.5){\makebox(2,1)[l]{4}}
\put(6,2){\circle*{0.7}}     \put(6.7,1.5){\makebox(2,1)[l]{6}}
\put(8,14){\circle*{0.7}}    \put(7.6,15){\makebox(2,1)[l]{2}}

\put(4.5,11.5){\makebox(2,1)[l]{\bf x}}
\put(4.5,7.5){\makebox(2,1)[l]{\bf z}}
\put(1.5,4){\makebox(2,1)[l]{\bf x}}
\put(5.8,3.5){\makebox(2,1)[l]{\bf y}}

\end{picture}
    \vline
\begin{picture}(10,18)
\put(4,14){\line(0,-1){8}}
\put(4,6){\line(-1,-2){2}}
\put(4,6){\line(1,-2){2}}
\put(4,14){\circle*{0.7}}    \put(3.5,15){\makebox(2,1)[l]{1}}
\put(4,10){\circle*{0.7}}    \put(4.7,9.5){\makebox(2,1)[l]{5}}
\put(4,6){\circle*{0.7}}     \put(4.7,5.5){\makebox(2,1)[l]{3}}
\put(2,2){\circle*{0.7}}     \put(2.7,1.5){\makebox(2,1)[l]{4}}
\put(6,2){\circle*{0.7}}     \put(6.7,1.5){\makebox(2,1)[l]{6}}
\put(8,14){\circle*{0.7}}    \put(7.6,15){\makebox(2,1)[l]{2}}

\put(4.5,11.5){\makebox(2,1)[l]{\bf x}}
\put(4.5,7.5){\makebox(2,1)[l]{\bf y}}
\put(1.5,4){\makebox(2,1)[l]{\bf x}}
\put(5.8,3.5){\makebox(2,1)[l]{\bf z}}
\end{picture}
    \vline
\begin{picture}(10,18)
\put(4,14){\line(-1,-2){2}}
\put(4,14){\line(1,-2){2}}
\put(8,14){\line(0,-1){8}}
\put(4,14){\circle*{0.7}}    \put(3.5,15){\makebox(2,1)[l]{1}}
\put(2,10){\circle*{0.7}}    \put(1.6,8){\makebox(2,1)[l]{4}}
\put(6,10){\circle*{0.7}}    \put(5.6,8){\makebox(2,1)[l]{6}}
\put(8,14){\circle*{0.7}}    \put(7.6,15){\makebox(2,1)[l]{2}}
\put(8,10){\circle*{0.7}}    \put(8.7,9.5){\makebox(2,1)[l]{3}}
\put(8,6){\circle*{0.7}}     \put(7.6,4){\makebox(2,1)[l]{5}}
\put(1.5,12){\makebox(2,1)[l]{\bf x}}
\put(5.5,12){\makebox(2,1)[l]{\bf y}}
\put(8.5,11.5){\makebox(2,1)[l]{\bf y}}
\put(8.5,7.5){\makebox(2,1)[l]{\bf x}}
\end{picture}
   \vline
\begin{picture}(10,18)
\put(4,14){\line(-1,-2){2}}
\put(4,14){\line(1,-2){2}}
\put(8,14){\line(0,-1){8}}
\put(4,14){\circle*{0.7}}    \put(3.5,15){\makebox(2,1)[l]{1}}
\put(2,10){\circle*{0.7}}    \put(1.6,8){\makebox(2,1)[l]{4}}
\put(6,10){\circle*{0.7}}    \put(5.6,8){\makebox(2,1)[l]{6}}
\put(8,14){\circle*{0.7}}    \put(7.6,15){\makebox(2,1)[l]{2}}
\put(8,10){\circle*{0.7}}    \put(8.7,9.5){\makebox(2,1)[l]{3}}
\put(8,6){\circle*{0.7}}     \put(7.6,4){\makebox(2,1)[l]{5}}
\put(1.5,12){\makebox(2,1)[l]{\bf x}}
\put(5.5,12){\makebox(2,1)[l]{\bf y}}
\put(8.5,11.5){\makebox(2,1)[l]{\bf x}}
\put(8.5,7.5){\makebox(2,1)[l]{\bf z}}
\end{picture}
   \vline
\begin{picture}(10,18)
\end{picture}
\\ \hline
\begin{picture}(11,11)
\put(3,7){\line(-1,-2){2}}
\put(3,7){\line(1,-2){2}}
\put(11,7){\line(0,-1){4}}
\put(3,7){\circle*{0.7}}    \put(2.5,8){\makebox(2,1)[l]{1}}
\put(1,3){\circle*{0.7}}    \put(0.6,1){\makebox(2,1)[l]{4}}
\put(5,3){\circle*{0.7}}    \put(4.6,1){\makebox(2,1)[l]{6}}
\put(7,7){\circle*{0.7}}    \put(6.6,8){\makebox(2,1)[l]{2}}
\put(11,7){\circle*{0.7}}   \put(10.6,8){\makebox(2,1)[l]{3}}
\put(11,3){\circle*{0.7}}   \put(10.6,1){\makebox(2,1)[l]{5}}
\put(0.5,5){\makebox(2,1)[l]{\bf x}}
\put(4.5,5){\makebox(2,1)[l]{\bf y}}
\put(11.5,4.5){\makebox(2,1)[l]{\bf x}}
\end{picture}
\\ \hline
\end{tabular}
\caption{An example for the recurrence relation \eqref{eq-special}.}
\label{fig-color}
\end{center}
\end{figure}

Figure \ref{fig-color} gives an example for the case $n=6,\,r=3,\,k=3$. In order to
make the figure clear, we replace the colors set $\{\bf 1,\bf 2,\bf 3\}$ by
$\{\bf x,\bf y,\bf z\}$.\vskip 0.3cm
\begin{cor}\label{cor-special}
The number of forests of $r$ special $k$-edge colored trees on $[n]$ with roots
$1,\ldots,r$ is
\begin{eqnarray*}
&&r(k-1)(kn-2n-r+1)(kn-2n-r+2)\cdots(kn-n-1)\\[5pt]
&&=r(k-1)(n-r-1)!{kn-n-1\choose n-r-1}.
\end{eqnarray*}
\end{cor}

\noindent{\it Proof.} It follows immediately from \eqref{eq-special} and the fact
$|\mathcal{E}_{n,k}^{n-1}[n\prec 1]|=k-1$.\qed

The following result is due to Gessel (see \cite[pp. 124]{Stanley}), and can also be deduced
in the same way as before.


\begin{cor}\label{cor-Gessel}
The number of $k$-edge colored trees on $[n]$ is
\[
k(n-2)!{kn-n\choose n-2}.
\]
\end{cor}

\noindent{\it Proof.} The case $r=1$ of Corollary \ref{cor-special} implies that the number of
special $k$-colored trees with root $1$ is
\[
(k-1)(n-2)!{kn-n-1\choose n-2}.
\]
Note that a $k$-colored tree is not special if and only if
one of the edges out of the root is colored by ${\bf k}$, and we can delete this edge to obtain a forest of two special $k$-colored trees such that one of the roots is $1$. But Corollary
\ref{cor-special} says that the number of forests of two special $k$-edge colored trees
with roots $1,2$ is
\[
2(k-1)(n-3)!{kn-n-1\choose n-3},
\]
so is the number of forests with roots $1,x$, where $x\neq 1$. Hence, the number of
$k$-colored trees on $[n]$ (with root $1$) is given by
\begin{align*}
\nonumber
&\hskip -2mm
(k-1)(n-2)!{kn-n-1\choose n-2}+2(n-1)(k-1)(n-3)!{kn-n-1\choose n-3}\\[5pt]
&=k(kn-n)(kn-n-1)\cdots(kn-2n+3).
\end{align*}
This completes the proof.\qed

\begin{cor}\label{cor-cor-deg}
The number of $k$-edge colored trees on $[n]$ with a specific root and root degree
$r$ is
\[
k(n-2)!{k-1\choose r-1}{(k-1)(n-1)\choose n-r-1}.
\]
\end{cor}

\noindent{\it Proof.} Suppose that the root of the tree is $1$. Then there are ${n-1\choose k}$
ways to choose the children of $1$ from the set $\{2,\ldots,n\}$, and
${k\choose r}r!$ ways to color the edges out of the root.
Obviously, each of subtrees of $1$ may be regarded as a special $k$-edge colored tree, because the
edges out of its root are colored from a definite set of $k-1$ colors. It follows
from Corollary \ref{cor-special} that there are $r(k-1)(n-r-2)!{k(n-1)-n\choose n-r-2}$
ways to construct the forest of $r$ special $k$-edge trees on $n-1$ vertices with $r$
specific roots. Therefore, the number of the required trees is
\[
{n-1\choose r}{k\choose r}r!r(k-1)(n-r-2)!{k(n-1)-n\choose n-r-2},
\]
as desired.\qed

\vskip 5mm \noindent{\bf Acknowledgments.} The authors would like to thank the anonymous referees for
helpful comments on a previous version of this paper. The first author was partially supported by the Natural Science Research Project of Ordinary
Universities in Jiangsu Province of China (grant no. 13KJB110001).


\vskip 1cm
\begin{thebibliography}{99}
\small \setlength{\itemsep}{-.8mm}

\bibitem{AZ}M. Aigner and G.M. Ziegler, Proofs from The Book, Fourth Ed., Springer-Verlag, Berlin, 2010.

\bibitem{Austin}
T.L. Austin, The enumeration of point labelled chromatic graphs and trees, Canad. J. Math. 12 (1960), 535--545.

\bibitem{Cayley}
A. Cayley, A theorem on trees, Quart. J. Math. 23 (1889), 376--378.

\bibitem{Chen90}
W.Y.C. Chen, A general bijective algorithm for trees, Proc. Natl. Acad. Sci. USA 87 (1990), 9635--9639.

\bibitem{ChenGuo}
W.Y.C. Chen and V.J.W. Guo, Bijections behind the Ramanujan polynomials, Adv. Appl. Math. (2001), 336--356.

\bibitem{ChenPeng}W.Y.C. Chen and J.F.F Peng, Disposition polynomials and plane trees, European J. Combin. 36 (2014), 122--129.

\bibitem{Clarke}
L.E. Clarke, On Cayley's formula for counting trees, J. London Math. Soc. 33 (1958), 26--28.

\bibitem{DuYin}R.R.X. Du and J. Yin, Counting labelled trees with given indegree sequence, J. Combin. Theory Ser. A 117 (2010), 345--353.

\bibitem{ER86}O. E\v{g}ecio\v{g}lu and J.B. Remmel, Bijections for Cayley trees,
spanning trees, and their $q$-analogues, J. Combin. Theory Ser. A 42 (1986), 15--30.

\bibitem{ER94}O. E\v{g}ecio\v{g}lu and J.B. Remmel, A bijection  for spanning trees of
complete multipartite graphs, Congr. Numer. 100 (1994), 225--243.

\bibitem{EE}
A. Erd$\acute{\rm e}$lyi and I.M.H. Etherington, Some problems of non-associative combinations II, Edinburgh Math. Notes 32 (1941), 7--12.

\bibitem{Fiedler}
M. Fiedler and J. Sedl$\acute{\rm a}\check{\rm c}$ek, $\check{\rm C}$asopis pro P$\check{\rm e}{\rm stov}\acute{\rm a}$ni Matematiky
83 (1958), 214--225.

\bibitem{GZ}V.J.W. Guo and J. Zeng, A generalization of the Ramanujan polynomials and plane trees, Adv. Appl. Math.  39 (2007), 96--115.

\bibitem{HS} M. Haiman and W. Schmitt, Incidence algebra antipodes and Lagrange inversion
in one and several variables, J. Combin. Theory Ser. A 50 (1989), 172--185.

\bibitem{Hou}Q.-H. Hou, An insertion algorithm and leaders of rooted trees, European J. Combin.  53  (2016), 35--44.

\bibitem{Narayana}T.V. Narayana, A partial order and its application to probability, Sankhy\={a} 21 (1959), 91--98.

\bibitem{Pitman}J. Pitman, Coalescent random forests, J. Combin. Theory, Ser. A 85 (1999), 165--193.

\bibitem{Renyi}A. R\'enyi, Some remarks on the theory of trees, MTA Mat. Kut. Inst. Kozl.
(Publ. math. Inst. Hungar. Acad. Sci.) 4 (1959), 73-85; Selected Papers Vol. 2,
Akad¨¦miai Kiad¨®, Budapest 1976, 363--374.

\bibitem{Riordan}J. Riordan, Forests of labeled trees, J. Combin. Theory 5 (1968), 90--103.

\bibitem{SZ}H. Shin and J. Zeng, A bijective enumeration of labeled trees with given indegree sequence,
J. Combin. Theory Ser. A 118 (2011), 115--128.


\bibitem{Stanley}R.P. Stanley, Enumerative Combinatorics Vol. 2, Cambridge University Press, Cambridge, 1999.
\end{thebibliography}
\end{document}